\documentclass{amsart}
\usepackage{amsmath}
\usepackage{amssymb}
\usepackage{amsfonts}

\newtheorem{theorem}{Theorem}
\theoremstyle{plain}

\newtheorem{definition}{Definition}

\newtheorem{remark}{Remark}

\numberwithin{equation}{section}
\input{tcilatex}

\begin{document}
\title[Univariate Right Caputo Fractional Ostrowski Inequalities]{Ostrowski
type inequalities involving the right Caputo fractional derivatives belong
to $L_{p}$}
\author{Mehmet Zeki SARIKAYA}
\email{sarikayamz@gmail.com}

\begin{abstract}
In this paper, we have established Ostrowski type inequalities involving the
right Caputo fractional derivatives belong to $L_{p}$ spaces ($1\leq p\leq
\infty $) via the right Caputo fractional Taylor formula with integral
remainder.
\end{abstract}

\maketitle

\section{Introduction}

In 1938, the classical integral \ inequality proved by A.M. Ostrowski as the
following:

\begin{theorem}
\label{t1} Let $f:\left[ a,b\right] \rightarrow 
\mathbb{R}
$ be continuous on $\left[ a,b\right] $ and differentiable on $\left(
a,b\right) $ whose derivative $f^{\prime }:\left( a,b\right) \rightarrow 
\mathbb{R}
$ is bounded on $\left( a,b\right) ,$ i.e., $\left\Vert f^{\prime
}\right\Vert _{\infty }=\underset{t\in \left( a,b\right) }{sup}\left\vert
f^{\prime }\left( t\right) \right\vert <+\infty .$ Then the following
inequality holds:%
\begin{equation}
\left\vert f(x)-\frac{1}{b-a}\int_{a}^{b}f(u)du\right\vert \leq \left(
b-a\right) \left[ \frac{1}{4}+\frac{\left( x-\frac{a+b}{2}\right) ^{2}}{%
\left( b-a\right) ^{2}}\right] \left\Vert f^{\prime }\right\Vert _{\infty }
\label{e.1.1}
\end{equation}%
for any $x\in \left[ a,b\right] $. The constant $\frac{1}{4}$ is the best
possible.
\end{theorem}

In \cite{Anastassiou1}, \ Anastassiou established a new Ostrowski inequality
which holds higer order derivatives functions.

\begin{theorem}
\label{tt1} Let $f\in C^{n+1}\left( \left[ a,b\right] \right) ,\ n\in 
\mathbb{N}$ and $x\in \left[ a,b\right] $ be fixed, such that $f^{(k)}(0)=0$%
, $k=1,...,n$. Then it holds 
\begin{equation}
\left\vert f(x)-\frac{1}{b-a}\int_{a}^{b}f(u)du\right\vert \leq \frac{%
\left\Vert f^{(n+1)}\right\Vert _{\infty }}{\left( n+2\right) !}\left[ \frac{%
\left( x-a\right) ^{n+2}+\left( b-x\right) ^{n+2}}{\left( b-a\right) }\right]
.  \label{e.1.3}
\end{equation}
\end{theorem}

Clearly inequality (\ref{e.1.3}) generalizes inequality (\ref{e.1.1}) for
higher order derivatives of $f$.

We give some necessary definitions and mathematical preliminaries of
fractional calculus theory \ which are used throughout this paper.

\begin{definition}
Let $f\in L_{1}[a,b].$ The Riemann-Liouville integrals $J_{a+}^{\alpha }f$
and $J_{b-}^{\alpha }f$ of order $\alpha >0$ with $a\geq 0$ are defined by 
\begin{equation*}
J_{a+}^{\alpha }f(x)=\frac{1}{\Gamma (\alpha )}\int_{a}^{x}\left( x-t\right)
^{\alpha -1}f(t)dt,\ \ x>a
\end{equation*}%
and%
\begin{equation*}
J_{b-}^{\alpha }f(x)=\frac{1}{\Gamma (\alpha )}\int_{x}^{b}\left( t-x\right)
^{\alpha -1}f(t)dt,\ \ x<b
\end{equation*}%
respectively where $\Gamma (\alpha )=\int_{0}^{\infty }e^{-t}u^{\alpha -1}du$%
. Here is $J_{a+}^{0}f(x)=J_{b-}^{0}f(x)=f(x).$
\end{definition}

\begin{definition}
Let $f\in AC^{m}\left( [a,b]\right) $($f^{(m-1)}$ is in $AC\left(
[a,b]\right) $), $\ m\in \mathbb{N}$, $m=[\alpha ],\ \alpha >0$ ($[.]$ the
ceiling of the number). The right Caputo fractional derivative of order $%
\alpha >0$ is defined by%
\begin{equation*}
D_{b-}^{\alpha }f(x)=\frac{(-1)^{m}}{\Gamma (m-\alpha )}\int_{x}^{b}\left(
t-x\right) ^{m-\alpha -1}f^{(m)}(t)dt,\ \ x\leq b.
\end{equation*}%
If $\alpha =m\in \mathbb{N}$, then%
\begin{equation*}
D_{b-}^{m}f(x)=(-1)^{m}f^{(m)}(x),\ \ \forall x\in \left[ a,b\right] .
\end{equation*}%
If $x>b,$ we define $D_{b-}^{\alpha }f(x)=0.$
\end{definition}

Properties concerning this operator can be found (\cite{Gorenflo}-\cite%
{samko}). For some recent results connected with fractional integral
inequalities see (\cite{Anastassiou1}-\cite{samko}).

In order to prove our main results, we need the following theorem proved by
Anastassiou in \cite{Anastassiou3}.

\begin{theorem}
\label{L1} Let $f\in AC^{m}\left( [a,b]\right) ,\ x\in \left[ a,b\right] $, $%
m=[\alpha ],\ \alpha >0.$ Then%
\begin{equation*}
f(x)=\dsum\limits_{k=0}^{m-1}\frac{f^{(k)}(b)}{k!}\left( x-b\right) ^{k}+%
\frac{1}{\Gamma (\alpha )}\int_{x}^{b}\left( t-x\right) ^{m-\alpha
-1}D_{b-}^{m}f(t)dt
\end{equation*}%
the right Caputo fractional Taylor formula with integral remainder.
\end{theorem}

In \cite{Anastassiou5}, \ Anastassiou established general univariate right
Caputo fractional Ostrowski inequalities with respect to $\left\Vert
.\right\Vert _{p},\ 1\leq p\leq \infty .$

\begin{theorem}
Let $f\in AC^{m}\left( [a,b]\right) ,$ $m=[\alpha ]$ and $f^{(k)}(b)=0,\
k=1,...,m-1.$ Then%
\begin{eqnarray}
&&\left\vert f(b)-\frac{1}{b-a}\int_{a}^{b}f(u)du\right\vert  \label{z} \\
&&  \notag \\
&\leq &\left\{ 
\begin{array}{cc}
\frac{\left\Vert D_{b-}^{m}f\right\Vert _{\infty ,\left[ a,b\right] }}{%
\Gamma (\alpha +2)}(b-a)^{\alpha }, & \text{if }D_{b-}^{\alpha }f\in
L_{\infty }\left( \left[ a,b\right] \right) ,\ \alpha >0 \\ 
\frac{\left\Vert D_{b-}^{m}f\right\Vert _{L_{1}\left( \left[ a,b\right]
\right) }}{\Gamma (\alpha +1)}(b-a)^{\alpha -1} & \text{if }D_{b-}^{\alpha
}f\in L_{1}\left( \left[ a,b\right] \right) ,~\alpha \geq 1 \\ 
\frac{\left\Vert D_{b-}^{m}f\right\Vert _{L_{q}\left( \left[ a,b\right]
\right) }}{\Gamma (\alpha )\left( p(\alpha -1)+1\right) ^{\frac{1}{p}%
}(\alpha +\frac{1}{p})}(b-a)^{\alpha -1+\frac{1}{p}} & 
\begin{array}{c}
\text{if }D_{b-}^{\alpha }f\in L_{q}\left( \left[ a,b\right] \right) , \\ 
\ p,q>1,\ \frac{1}{p}+\frac{1}{q}=1,\ \alpha >1-\frac{1}{p}%
\end{array}%
\end{array}%
\right.  \notag
\end{eqnarray}
\end{theorem}

The aim of this paper is to establish Ostrowski type inequalities involving
the right Caputo fractional derivatives belong to $L_{p}$ spaces ($1\leq
p\leq \infty $) via the right Caputo fractional Taylor formula with integral
remainder.

\section{Main Results}

We start with the following theorem:

\begin{theorem}
\label{T1} Let $f,g\in AC^{m}\left( [a,b]\right) ,$ $m=[\alpha ],\ \alpha >0$%
. Assume $f^{(k)}(b)=g^{(k)}(b)=0,\ k=1,...,m-1$ and $D_{b-}^{\alpha
}f,D_{b-}^{\alpha }g\in L_{\infty }\left( \left[ a,b\right] \right) .$ Then%
\begin{eqnarray}
&&\left\vert 2\int_{a}^{b}f(x)g(x)dx-\int_{a}^{b}\left(
f(x)g(b)+g(x)f(b)\right) dx\right\vert  \label{A} \\
&&  \notag \\
&\leq &\left\Vert D_{b-}^{\alpha }f\right\Vert _{\infty }J_{a+}^{\alpha
+1}\left\vert g(b)\right\vert +\left\Vert D_{b-}^{\alpha }g\right\Vert
_{\infty }J_{a+}^{\alpha +1}\left\vert f(b)\right\vert .  \notag
\end{eqnarray}
\end{theorem}

\begin{proof}
Let $x\in \lbrack a,b].$ Using Theorem \ref{L1} and from the hypothesis of
Theorem \ref{T1}, we have the following identities 
\begin{equation}
f(x)-f(b)=\frac{1}{\Gamma (\alpha )}\int_{x}^{b}\left( t-x\right) ^{\alpha
-1}D_{b-}^{\alpha }f(t)dt  \label{E4}
\end{equation}%
\begin{equation}
g(x)-g(b)=\frac{1}{\Gamma (\alpha )}\int_{x}^{b}\left( t-x\right) ^{\alpha
-1}D_{b-}^{\alpha }g(t)dt.  \label{E5}
\end{equation}%
Multiplying both sides of (\ref{E4}) and (\ref{E5}) by $g(x)$ and $f(x)$
respectively and adding the resulting identities, we have%
\begin{eqnarray*}
&&2f(x)g(x)-f(b)g(x)-f(x)g(b) \\
&=&\frac{g(x)}{\Gamma (\alpha )}\int_{x}^{b}\left( t-x\right) ^{\alpha
-1}D_{b-}^{\alpha }f(t)dt+\frac{f(x)}{\Gamma (\alpha )}\int_{x}^{b}\left(
t-x\right) ^{\alpha -1}D_{b-}^{\alpha }g(t)dt.
\end{eqnarray*}%
Integrating the resulting inequality with respest to $x$ over $[a,b]$ and
using the properties of modulus, we obtain%
\begin{eqnarray}
&&\left\vert 2\int_{a}^{b}f(x)g(x)dx-\int_{a}^{b}\left(
f(b)g(x)+f(x)g(b)\right) dx\right\vert  \notag \\
&&  \label{E} \\
&\leq &\int_{a}^{b}\frac{\left\vert g(x)\right\vert }{\Gamma (\alpha )}%
\left( \int_{x}^{b}\left( t-x\right) ^{\alpha -1}\left\vert D_{b-}^{\alpha
}f(t)\right\vert dt\right) dx+\int_{a}^{b}\frac{\left\vert f(x)\right\vert }{%
\Gamma (\alpha )}\left( \int_{x}^{b}\left( t-x\right) ^{\alpha -1}\left\vert
D_{b-}^{\alpha }g(t)\right\vert dt\right) dx.  \notag
\end{eqnarray}%
Hence it holds%
\begin{eqnarray*}
&&\left\vert 2\int_{a}^{b}f(x)g(x)dx-\int_{a}^{b}\left(
f(b)g(x)+f(x)g(b)\right) dx\right\vert \\
&& \\
&\leq &\frac{\left\Vert D_{b-}^{\alpha }f\right\Vert _{\infty }}{\Gamma
(\alpha )}\int_{a}^{b}\left\vert g(x)\right\vert \left( \int_{x}^{b}\left(
t-x\right) ^{\alpha -1}dt\right) dx+\frac{\left\Vert D_{b-}^{\alpha
}g\right\Vert _{\infty }}{\Gamma (\alpha )}\int_{a}^{b}\left\vert
f(x)\right\vert \left( \int_{x}^{b}\left( t-x\right) ^{\alpha -1}dt\right) dx
\\
&& \\
&=&\frac{\left\Vert D_{b-}^{\alpha }f\right\Vert _{\infty }}{\Gamma (\alpha
+1)}\int_{a}^{b}\left( b-x\right) ^{\alpha }\left\vert g(x)\right\vert dx+%
\frac{\left\Vert D_{b-}^{\alpha }g\right\Vert _{\infty }}{\Gamma (\alpha +1)}%
\int_{a}^{b}\left( b-x\right) ^{\alpha }\left\vert f(x)\right\vert dx \\
&& \\
&=&\left\Vert D_{b-}^{\alpha }f\right\Vert _{\infty }J_{a+}^{\alpha
+1}\left\vert g(b)\right\vert +\left\Vert D_{b-}^{\alpha }g\right\Vert
_{\infty }J_{a+}^{\alpha +1}\left\vert f(b)\right\vert
\end{eqnarray*}%
which the proof is completed.
\end{proof}

\begin{remark}
If we take $g(x)=1$ in Theorem \ref{T1}, the inequality (\ref{A}) reduces
the first inequality in (\ref{z}).
\end{remark}

\begin{theorem}
\label{T2} Let $f,g\in AC^{m}\left( [a,b]\right) ,$ $m=[\alpha ],\ \alpha
\geq 1$. Assume $f^{(k)}(b)=g^{(k)}(b)=0,\ k=1,...,m-1$ and $D_{b-}^{\alpha
}f,D_{b-}^{\alpha }g\in L_{1}\left( \left[ a,b\right] \right) .$ Then%
\begin{eqnarray}
&&\left\vert 2\int_{a}^{b}f(x)g(x)dx-\int_{a}^{b}\left(
f(x)g(b)+g(x)f(b)\right) dx\right\vert  \label{A1} \\
&&  \notag \\
&\leq &\left\Vert D_{b-}^{\alpha }f\right\Vert _{L_{1}\left( \left[ a,b%
\right] \right) }J_{a+}^{\alpha }\left\vert g(b)\right\vert +\left\Vert
D_{b-}^{\alpha }g\right\Vert _{L_{1}\left( \left[ a,b\right] \right)
}J_{a+}^{\alpha }\left\vert f(b)\right\vert .  \notag
\end{eqnarray}
\end{theorem}

\begin{proof}
From the inequality (\ref{E}) of Theorem \ref{T1}, we have again%
\begin{eqnarray*}
&&\left\vert 2\int_{a}^{b}f(x)g(x)dx-\int_{a}^{b}\left(
f(b)g(x)+f(x)g(b)\right) dx\right\vert \\
&& \\
&\leq &\int_{a}^{b}\frac{\left\vert g(x)\right\vert }{\Gamma (\alpha )}%
\left( \int_{x}^{b}\left( t-x\right) ^{\alpha -1}\left\vert D_{b-}^{\alpha
}f(t)\right\vert dt\right) dx+\int_{a}^{b}\frac{\left\vert f(x)\right\vert }{%
\Gamma (\alpha )}\left( \int_{x}^{b}\left( t-x\right) ^{\alpha -1}\left\vert
D_{b-}^{\alpha }g(t)\right\vert dt\right) dx \\
&& \\
&\leq &\int_{a}^{b}\frac{\left\vert g(x)\right\vert }{\Gamma (\alpha )}%
\left( b-x\right) ^{\alpha -1}\left( \int_{x}^{b}\left\vert D_{b-}^{\alpha
}f(t)\right\vert dt\right) dx+\int_{a}^{b}\frac{\left\vert f(x)\right\vert }{%
\Gamma (\alpha )}\left( b-x\right) ^{\alpha -1}\left( \int_{x}^{b}\left\vert
D_{b-}^{\alpha }g(t)\right\vert dt\right) dx \\
&& \\
&=&\left\Vert D_{b-}^{\alpha }f\right\Vert _{L_{1}\left( \left[ a,b\right]
\right) }\int_{a}^{b}\frac{\left\vert g(x)\right\vert }{\Gamma (\alpha )}%
\left( b-x\right) ^{\alpha -1}dx+\left\Vert D_{b-}^{\alpha }g\right\Vert
_{L_{1}\left( \left[ a,b\right] \right) }\int_{a}^{b}\frac{\left\vert
f(x)\right\vert }{\Gamma (\alpha )}\left( b-x\right) ^{\alpha -1}dx \\
&& \\
&=&\left\Vert D_{b-}^{\alpha }f\right\Vert _{L_{1}\left( \left[ a,b\right]
\right) }J_{a+}^{\alpha }\left\vert g(b)\right\vert +\left\Vert
D_{b-}^{\alpha }g\right\Vert _{L_{1}\left( \left[ a,b\right] \right)
}J_{a+}^{\alpha }\left\vert f(b)\right\vert
\end{eqnarray*}%
which completes the proof.
\end{proof}

\begin{remark}
If we take $g(x)=1$ in Theorem \ref{T2}, the inequality (\ref{A1}) reduces
the second inequality in (\ref{z}).
\end{remark}

\begin{theorem}
\label{h1} Let $f,g\in AC^{m}\left( [a,b]\right) ,$ $m=[\alpha ],\ \alpha
>0,\ p,q>1,\ \frac{1}{p}+\frac{1}{q}=1,\ \alpha >1-\frac{1}{p}$. Assume $%
f^{(k)}(b)=g^{(k)}(b)=0,\ k=1,...,m-1$ and $D_{b-}^{\alpha }f,D_{b-}^{\alpha
}g\in L_{q}\left( \left[ a,b\right] \right) .$ Then%
\begin{eqnarray}
&&\left\vert 2\int_{a}^{b}f(x)g(x)dx-\int_{a}^{b}\left(
f(x)g(b)+g(x)f(b)\right) dx\right\vert  \notag \\
&&  \label{A2} \\
&\leq &\Gamma \left( \alpha +\frac{1}{p}\right) \left( \left\Vert
D_{b-}^{\alpha }f\right\Vert _{L_{q}\left( \left[ a,b\right] \right)
}J_{a+}^{\alpha +\frac{1}{p}}\left\vert g(b)\right\vert +\left\Vert
D_{b-}^{\alpha }g\right\Vert _{L_{q}\left( \left[ a,b\right] \right)
}J_{a+}^{\alpha +\frac{1}{p}}\left\vert f(b)\right\vert \right) .  \notag
\end{eqnarray}
\end{theorem}

\begin{proof}
From the inequality (\ref{E}) of Theorem \ref{T1}, we have again%
\begin{eqnarray*}
&&\left\vert 2\int_{a}^{b}f(x)g(x)dx-\int_{a}^{b}\left(
f(b)g(x)+f(x)g(b)\right) dx\right\vert \\
&& \\
&\leq &\int_{a}^{b}\frac{\left\vert g(x)\right\vert }{\Gamma (\alpha )}%
\left( \int_{x}^{b}\left( t-x\right) ^{\alpha -1}\left\vert D_{b-}^{\alpha
}f(t)\right\vert dt\right) dx+\int_{a}^{b}\frac{\left\vert f(x)\right\vert }{%
\Gamma (\alpha )}\left( \int_{x}^{b}\left( t-x\right) ^{\alpha -1}\left\vert
D_{b-}^{\alpha }g(t)\right\vert dt\right) dx.
\end{eqnarray*}%
Using the H\"{o}lder's integral inequality, we obtain%
\begin{eqnarray*}
&&\left\vert 2\int_{a}^{b}f(x)g(x)dx-\int_{a}^{b}\left(
f(b)g(x)+f(x)g(b)\right) dx\right\vert \\
&& \\
&\leq &\int_{a}^{b}\frac{\left\vert g(x)\right\vert }{\Gamma (\alpha )}%
\left( \int_{x}^{b}\left( t-x\right) ^{(\alpha -1)p}dt\right) ^{\frac{1}{p}%
}\left( \int_{x}^{b}\left\vert D_{b-}^{\alpha }f(t)\right\vert ^{q}dt\right)
^{\frac{1}{q}}dx \\
&& \\
&&+\int_{a}^{b}\frac{\left\vert f(x)\right\vert }{\Gamma (\alpha )}\left(
\int_{x}^{b}\left( t-x\right) ^{(\alpha -1)p}dt\right) ^{\frac{1}{p}}\left(
\int_{x}^{b}\left\vert D_{b-}^{\alpha }g(t)\right\vert ^{q}dt\right) ^{\frac{%
1}{q}}dx \\
&& \\
&=&\frac{\left\Vert D_{b-}^{\alpha }f\right\Vert _{L_{q}\left( \left[ a,b%
\right] \right) }}{\Gamma (\alpha )}\int_{a}^{b}\left\vert g(x)\right\vert
\left( \int_{x}^{b}\left( t-x\right) ^{(\alpha -1)p}dt\right) ^{\frac{1}{p}%
}dx \\
&& \\
&&+\frac{\left\Vert D_{b-}^{\alpha }g\right\Vert _{L_{q}\left( \left[ a,b%
\right] \right) }}{\Gamma (\alpha )}\int_{a}^{b}\left\vert f(x)\right\vert
\left( \int_{x}^{b}\left( t-x\right) ^{(\alpha -1)p}dt\right) ^{\frac{1}{p}}
\end{eqnarray*}%
By simple computation%
\begin{equation*}
\left( \int_{x}^{b}\left( t-x\right) ^{(\alpha -1)p}dt\right) ^{\frac{1}{p}}=%
\frac{(b-x)^{\alpha -1+\frac{1}{p}}}{\left( p(\alpha -1)+1\right) ^{\frac{1}{%
p}}}.
\end{equation*}%
This last equality is substituted the above, then we have the conclusion.
\end{proof}

\begin{remark}
If we take $g(x)=1$ in Theorem \ref{h1}, the inequality (\ref{A2}) reduces
the third inequality in (\ref{z}).
\end{remark}

\end{document}